\documentclass[12pt,a4paper]{article}
  
\usepackage{amsmath,amssymb} 
\usepackage{breakurl}
\usepackage[breaklinks]{hyperref}
\usepackage[small]{caption}
\usepackage{euscript,upref} 
\usepackage{cite} 
\usepackage{url} 
  
\newcommand{\mbf}[1]{\protect{\makebox{\boldmath$#1$}}}
\newcommand{\mbb}{\mathbb}

\newcommand{\eus}{\EuScript} 
\newcommand{\dist}{\mathrm{dist}\;} 
\newcommand{\ov}{\overline} 
\newcommand{\un}{\underline}

\newcommand{\dual}{\mathrm{dual}\,}

\title{\bf Non-traditional intervals and their use.\\[2mm] 
           Which ones really make sense?} 

\author{\sf Sergey P. Shary} 
\date{\small Federal Research Center for Information and Computational Technologies\\ 
and Novosibirsk State University,\, Novosibirsk,\, Russia\\[1mm] 
E-mail: \url{shary@ict.nsc.ru}} 
  
\begin{document}
\maketitle
  
\begin{abstract}
The paper discusses the question of why intervals, which are the main object of Interval 
Analysis, have exactly the form that we know well and habitually use, and not some other. 
In particular, we investigate why traditional intervals are closed, i.\,e. contain their 
endpoints, and also what is wrong with an empty interval. The second question considered 
in the work is how expedient it is to expand the set of traditional intervals by some 
other objects. We show that improper (``reversed'') intervals and the arithmetic of such 
intervals (Kaucher complete interval arithmetic) are very useful from many different 
points of view.  \\[2mm] 
\end{abstract}
  
\bigskip 
  
\section{Introduction} 
\label{IntroSect}

The article discusses one of the basic objects of Interval Analysis, namely, the concept 
of interval. Recall that classical intervals are closed, connected, and bounded subsets 
of the real line $\mbb{R}$, i.e., sets of the form 
\begin{equation} 
\label{UsualInterval}
[a, b] = \bigl\{\,x\in\mbb{R} \mid a\leq x\leq b\,\bigr\} 
\end{equation} 
(see \cite{AlefeldHerzberger, HansenWalster,MayerBook,MooreBakerCloud,NeumaierBook, 
SharyBook} and other books on Interval Analysis). The set of all intervals (usually 
with  arithmetic operations on it) is denoted by $\mbb{IR}$. Multidimensional intervals 
are their generalizations, in one sense or another. 
  
The issues covered below were raised in an online discussion on 
\verb|reliable_computing| \texttt{mailing list} (see \cite{RCmaillist}) that took 
place in the spring of 2018. Its starting point was the question of whether it is 
advisable to introduce and further use open and half-open intervals, such as $[a, b[\,$, 
$]a, b[\,$, $]a, b]$, in addition to the existing closed intervals \eqref{UsualInterval}  
(we denote various types of intervals in the style of N.\,Bourbaki). The author had 
experience of working with similar objects and therefore took an active part in 
the ensuing discussion. A summary of the views on these issues constitutes the core 
of this article. Previously, some ideas of the following text have been published 
in the book \cite{SharyBook}, as Section~1.11. 
  
Hereafter, by ``traditional intervals'', we mean usual intervals of the form 
\eqref{UsualInterval} that constitute the classical interval arithmetic $\mbb{IR}$, 
and ``non-traditional'' intervals in Section~\ref{ClosedIntsSect}  are  open  and 
half-open intervals. Further, in Section~\ref{OtherNonTraditSect},  we consider 
improper (``reversed'') intervals from Kaucher interval arithmetic $\mbb{KR}$ 
(algebraic and order completion of $\mbb{IR}$) as non-traditional intervals. 
A short Section~\ref{EmptyIntSect} discusses usefulness of a special ``empty 
interval''. Our notation follows the informal international standard \cite{INotation}. 
  
Infinite and semi-infinite intervals of the form $[-\infty, p]$, $[q, \infty]$, and 
$[-\infty, p]\cap[q, \infty]$ can also be classified as non-traditional. They have 
been first considered by W.\,Kahan \cite{Kahan-68} and further developed in detail 
in the works \cite{Laveuve, Ratz} and some others. Usually, the interval arithmetic 
of such infinite intervals is called ``Kahan interval arithmetic'', although other 
terms may also be used (see, e.\,g. \cite{KearfottBook}). Over the years, Kahan 
interval arithmetic has received numerous applications in Interval Analysis, and 
hence it does not need additional justification from our side. For this reason, 
we do not consider these non-traditional intervals in our work.

  
\section{Why, in Interval Analysis, intervals are closed?} 
\label{QuestionSect}

One of the popular myths about Interval Analysis and interval arithmetic $\mbb{IR}$ 
(widespread among beginners and veterans alike) is that it has the disadvantage of 
supporting only intervals which contain their endpoints. This allegedly makes it 
impossible to perform some basic set operations (like complement), bounds its 
expressive power, limits arithmetic operations, etc., etc. 
  
Well, is it really a drawback that intervals are closed and that they have no open 
endpoints?? Or, on the contrary, we do not understand the underlying reasons behind 
$\mbb{IR}$? 
  
This kind of questions was asked by W.\,Kahan more than half a century ago (see 
\cite{Kahan-68}), and he was the first to propose the introduction of non-closed 
intervals (although he did not realize his idea to the end). Then the same issues 
came up in 1998 and then in 2018. 
  
Of course, in the practical application of Interval Analysis, these questions are not 
so relevant. The fact is, all numerical values of continuously varying quantities 
encountered in practice are approximate, all measurements are performed with non-zero 
error, etc. Therefore, it is almost impossible to trace the difference between a point 
and other values that differ infinitesimally little from it. As a consequence, for 
engineering practice, the difference between open or closed intervals seems rather 
abstract and unimportant. But this question is important in theory and in calculations 
with intervals. 
  
There are arguments in favor of admitting non-traditional open and half-open intervals 
in computation and in our reasoning in general, in favor of giving them ``citizenship 
rights'' in Interval Analysis. Sometimes such intervals seem to significantly expand 
our capabilities. For example, let $[0, 1]$ be a traditional closed interval and let 
$]0, 1]$ be half-open at zero, then 
\begin{equation} 
\label{ZeroDivision} 
\frac{[0, 1]}{[0, 1]} = [-\infty, \infty],  \qquad
   \text{ but }\quad\frac{[0, 1]}{\,]0, 1]} = [0, \infty], 
\end{equation} 
since the first fraction must contain the result of the division $0/0$, and this is 
not the case for the second fraction. 
  
In other words, the benefit from the fact that, in \eqref{ZeroDivision}, the denominator 
is open at one of its endpoints is half the real line! This is very helpful in interval 
Newton method for enclosing zeros of functions. Its iterations are defined by the formula 
\begin{equation}
\label{INewtonIters} 
\mbf{X}^{(k+1)}\gets \mbf{X}^{(k)}\cap 
  {\,\eus{N}}\bigl(\mbf{X}^{(k)},\tilde{x}^{(k)}\bigr), 
  \quad  \tilde{x}^{(k)}\in\mbf{X}^{(k)}, 
  \qquad  k = 0,1,2, \ldots , 
\end{equation} 
where 
\begin{equation} 
\label{INewtonOperator} 
\rule[-6mm]{0mm}{14mm} 
{\eus{N}}(\mbf{X},\tilde{x}) 
  := \tilde{x} - \frac{f(\tilde{x})}{\mbf{f}'(\mbf{X})} 
\end{equation} 
is the interval Newton operator (see e.g. \cite{HansenWalster, MayerBook, 
MooreBakerCloud, NeumaierBook, SharyBook}). If, in the fraction  from 
\eqref{INewtonOperator}, the numerator and denominator coincide with those 
from \eqref{ZeroDivision}, then we would get very large improvement of the result. 
Then the interval Newton method finds solutions to equations much faster, areas 
that obviously did not contain solutions are better eliminated, etc. 
  
A similar situation can arise when implementing the Gauss-Seidel interval method 
(see e.g. \cite{HansenWalster, KearfottBook, MooreBakerCloud, NeumaierBook, SharyBook}). 
One also need to divide by an interval there and then intersect the result with another 
interval, and in this case we can again get the same huge improvement of the final 
result. 
  
More than 20 years ago, I worked for the Novosibirsk software company UniPro, which 
then carried out orders for Sun Microsystems, Inc.  Our team was implementing 
Sun's interval Fortran-95 (see \cite{SunIASpecification}), a programming language 
with built-in interval data type and operations with it. Just at that time, our project 
manager, Bill Walster from Sun Microsystems, was writing a book \cite{HansenWalster} 
with E.\,Hansen, which was devoted to interval methods for global optimization and 
equation solving. He was greatly impressed by the above observation with the interval 
Newton method. Hence, it was decided to implement something like open or half-open 
intervals, at least partially, insofar as the sign of the zero, which was inherent 
to floating-point computer arithmetics according to IEEE 754/854 standards, made it 
possible to easily implement open and closed endpoints at zero. 
  
In a computer, an interval $[a, b]$ is naturally represented by a pair of real numbers, 
$(a, b)$, and, for zero endpoints, we can take that 
\begin{center} 
\tabcolsep=5mm 
\begin{tabular}{lll} 
$(a, +0)\,$ means $\;[a, 0]$, & $(a, -0)\,$ means $\;[a, 0[$ &  for $a < 0$,\\[2mm] 
$(-0, b)\,$ means $\;[0, b]$, & $(+0, b)\,$ means $\;]0, b]$ &  for $b > 0$. 
\end{tabular} 
\end{center} 
Approximately the same was proposed for the implementation of intervals not closed 
at zero by W.\,Kahan in his work \cite{Kahan-68}. Recently, J.\,Gustafson proposed 
a modern computer implementation of open intervals based on the construction of 
\emph{universal numbers}, Unums, developed in \cite{Gustafson}. In general, 
the possibility of implementing open intervals makes available the modified division 
from \eqref{ZeroDivision} and all the bonuses that follow from it. As a result, 
in 1999, Bill Walster immediately jumped at the idea. 
  
In our team, I was responsible for semantic testing and general mathematical consulting 
of the works on interval Fortran-95, so I had to think deeply about the consequences 
of introducing a new construction into the language. As I delved into the question, 
I realized that implementation of open intervals, even partially, was hardly possible 
and did not make much sense. It turned out that there are very important and even 
fundamental mathematical reasons why the intervals should be closed. The fact is that, 
with their mathematical properties, bounded closed intervals essentially differ from 
non-closed, open and half-open, intervals. That was exactly what I reported to Bill 
Walster, and the implementation of non-closed intervals in Fortran-95 was canceled. 
  
Then the discussion in our team acquired a broader context, and the participants began 
to discuss questions about whether any other intervals, besides the classical intervals, 
are needed at all. Is it necessary to introduce the empty set into our algebraic systems 
of intervals? These topics are discussed in Sections~\ref{OtherNonTraditSect} and 
\ref{EmptyIntSect} of the work, but next we will look at what is so good about closed 
intervals.

  
\section{Closed intervals are compacts, \\* 
         complete lattices and complete metric spaces} 
\label{ClosedIntsSect}

Let us remind that a subset $S$ of a topological space $X$ is called \emph{compact}, 
if for every open cover of $S$ there exists a finite subcover of $S$ \cite{Dieudonne, 
Rudin,Shilov}. In fact, the concept of compactness formalizes the idea that it is 
possible to exhaust a set by a finite number of its arbitrarily small open subsets. 
Bounded closed intervals are compact sets in $\mbb{R}$, while open and half-open 
intervals are non-compact, with all the ensuing consequences. 
  
Considered as metric spaces, i.\,e. sets with an abstract distance (metric), non-closed 
intervals are not \emph{complete spaces} (see \cite{Dieudonne, Rudin}): the fundamental 
sequences (also called \emph{Cauchy sequences}) of elements from such non-closed 
intervals do not necessarily converge within such intervals. For example, the sequence 
of numbers $1/k$, $k = 1,2, \ldots\,$, has no limit within the half-open interval 
$\,]0, 1]$. 
  
Partially ordered sets in which we can freely take the supremum and infimum for each 
pair of elements is called \emph{lattice} \cite{Birkhoff}. A partially ordered set 
is called \emph{complete lattice} if all of its subsets have a supremum and an infimum. 
Non-closed intervals are not complete lattices with respect to the standard order 
``$\leq$'' on the real line, i.\,e. one cannot take infima and suprema, with respect 
to the order ``$\leq$'', for every subset of a non-closed interval. An example is 
the same sequence of numbers $1/k$, $k = 1,2, \ldots\,$, in the half-open interval 
$\,]0, 1]$. 
  
The above facts have many unpleasant consequences for practice, which we list below.

\paragraph{Non-compactness.}
On compact closed intervals, continuous functions reach their extrema, i.\,e. their 
minimal and maximal values (the Weierstrass extreme value theorem). But on non-compact, 
open and half-open intervals, continuous functions may not reach their extreme values.

\paragraph{Brouwer fixed point theorem and Banach fixed-point theorem.} 
The Brouwer's fixed point theorem (see \cite{GranasDugundji,OrtegaRheinboldt,Zeidler}) 
states that for any continuous function $\phi$  mapping a compact convex set  of 
$\mbb{R}^n$ to itself, there is such a point $\tilde{x}$ that $\tilde{x} = \phi(\tilde{x})$ 
(a ``fixed point'' that remains unchanged). Obtained in 1909--1912, this result has become 
one of the cornerstones of computational Interval Analysis, since intervals in $\mbb{R}$ 
and their multidimensional analogs are convex compact sets. Given an equation $f(x) = 0$, 
we can always reduce it to a fixed-point form $x = \phi(x)$ and then proceed as follows. 
If, using interval methods for enclosing the ranges of functions, we check the fulfillment 
of the Brouwer fixed-point theorem on an interval $\mbf{X}$ for the mapping $\phi$, i.\,e. 
that the inclusion $\phi(\mbf{X})\subseteq\mbf{X}$ is valid, then we rigorously prove 
the existence of a solution to the fixed-point equation $x = \phi(x)$ within $\mbf{X}$. 
  
The above technique, which constructively proves existence of solutions to equations, is 
an integral part of important interval methods, in particular, the Krawczyk method, 
the interval Newton method, and the Hansen-Sengupta method (see e.g. \cite{HansenWalster, 
KearfottBook, MayerBook, MooreBakerCloud, NeumaierBook, SharyBook}). 
  
The Banach fixed-point theorem (also known as the \emph{contraction mapping theorem}; 
see e.\,g. \cite{GranasDugundji, OrtegaRheinboldt, Zeidler}) states that a complete 
metric space $X$ with a contraction mapping $\phi: X\to X$ has a unique fixed-point 
$\tilde{x}$, i.\,e. such that $\tilde{x} = \phi(\tilde{x})$. It is also an important 
tool of computational Interval Analysis, because traditional intervals are complete 
metric spaces and hence the Banach fixed-point theorem allows one to prove 
the existence of solutions to equations and even their uniqueness. 
  
For non-closed intervals and their multidimensional analogs, the Brouwer fixed point 
theorem and the Banach fixed-point theorem are not valid. 
  
Therefore, the interval tests for the existence of solutions to equations and 
systems of equations, that is, the interval Newton method, the Krawczyk method, 
the Hansen-Sengupta method do not work in full with non-closed intervals. 
  
Thus, with non-closed intervals, computational Interval Analysis is deprived 
of its most powerful tools, widely used in solving equations, systems of linear and 
nonlinear equations as well as in global optimization.

\paragraph{The nested intervals principle.} 
It is known to be one of popular interval 
tools for both theory and verified computing. Let us recall its formulation: 
\begin{quote} 
Every nested interval sequence $\{\mbf{X}_{k}\}$, i.\,e. 
such that $\mbf{X}_{k+1}\subseteq \mbf{X}_{k}$      \\* 
for all $k$, converges and has the limit $\cap_{k=1}^{\infty} \mbf{X}_{k}$. 
\end{quote} 
$\mbf{X}_k$ may be either one-dimensional intervals or interval boxes in $\mbb{R}^n$. 
  
It turns out to be incorrect for non-traditional intervals. The sequence of half-open 
intervals does not necessarily converge to anything and may have an empty intersection. 
For example, 
\begin{equation*} 
\bigcap_{k=1}^\infty  \hspace{1.1ex} 
   \bigl]\hspace{0.1ex}0,\, \tfrac{1}{k}\,\bigr] \   = \  \varnothing. 
\end{equation*} 
  
This is a great loss. Let us recall that the theory of interval integral and 
interval estimates of the integral of a real function are based on this principle 
(see \cite{CapraniMadsenRall, Rall-82}). 
  
The most practical and efficient interval methods for solving operator equations 
(integral and differential) are based on the nested intervals principle and 
some fixed point theorems, in particular, the Banach fixed-point theorem and 
the Schr\"{o}der fixed-point theorem \cite{Collatz, NeumaierBook}. They also become 
invalid, since non-closed intervals are not complete topological spaces.

\paragraph{The Birkhoff-Tarski theorem and the Kantorovich lemma.} 
These are popular fixed-point theorems for partially ordered sets, analogs 
of topological fixed-point theorems that we formulated earlier. Birkhoff-Tarski 
principle (also called Knaster-Tarski principle; see \cite{Birkhoff, GranasDugundji, 
Tarski}) states that if $X$ is a complete lattice and $\phi: X\to X$ is an isotone 
(order preserving) function, then $\phi$ has a fixed point $\tilde{x}\in X$, i.\,e., 
such that $\tilde{x} = \phi(\tilde{x})$. A feature of this result is the absence of 
special requirements for the continuity of the function $\phi$, the form of the set
$X$, and its topological properties. 
  
Birkhoff-Tarski theorem is incorrect for non-closed intervals that are not complete 
lattices with respect to the standard order ``$\leq$'' on $\mbb{R}$ and with respect 
to inclusion ordering between intervals. 
  
Kantorovich lemma (see \cite{OrtegaRheinboldt}, Section 13.2) is a similar useful 
result, which is a fixed-point theorem for isotone mappings. It also becomes 
invalid for non-closed intervals.

\paragraph{Distance between various types of intervals.} 
How should we calculate the distance between $[\alpha, \beta[$ and $[\alpha, \beta]$, 
two intervals that differ in only one endpoint? 
  
In mathematics, distance (deviation) is usually formalized by the concept of 
a \emph{metric}, a function that gives a distance between each pair of elements of 
a set. The metric is defined axiomatically, as a nonnegative function that satisfies 
three axioms: identity of indiscernibles, symmetry, and triangle inequality 
(see details e.g. in \cite{Collatz, Dieudonne, Engelking}). 
  
The distance between intervals is known  to be determined as follows (see 
\cite{AlefeldHerzberger, MayerBook, MooreBakerCloud, NeumaierBook, SharyBook}) 
\begin{equation*} 
\dist(\mbf{a}, \mbf{b}) 
   = \max\bigl\{|\un{\mbf{a}} - \un{\mbf{b}}|, |\ov{\mbf{a}} - \ov{\mbf{b}}|\bigr\}. 
\end{equation*} 
It should be equal to zero for $\mbf{a} = [\alpha, \beta[$ and $\mbf{b} = 
[\alpha, \beta]$. Thus, one of the main purposes of distance, which is to distinguish 
between elements of a set that do not coincide with each other, is not fulfilled. 
Moreover, the metric (distance) can not be introduced in any way on the set of all 
closed and non-closed intervals , that is, this set is essentially non-metrizable 
as a topological space. 
  
The Arkhangelskii metrization theorem (see \cite{Engelking}) asserts that the topology 
of a space can be determined by a metric if and only if this space satisfies the first 
axiom of separation (the so-called T1-axiom) and has a countable fundamental family 
of open neighborhoods. Axiom T1 is the weakest axiom of separability, it requires that 
any one of two points of the space has a neighborhood not containing the other point. 
It is easy to see that the space of all closed and non-closed intervals does not satisfy 
even this weakest axiom: a half-open interval $[a, b[$ and its closure $[a, b]$ can not 
be surrounded by such neighborhoods. 
  
The failure of the axiom T1 is a very serious evidence of the fact that the topological 
space under consideration is very exotic, even pathologic. For us, in Interval Analysis, 
it implies that a meaningful calculus on the set of closed and non-closed intervals 
is most likely never to be constructed. Of course, this does not exclude individual 
episodic applications of non-closed intervals in certain particular situations. 
But in general, alas \ldots

  
\section{Can other non-traditional intervals \\* be useful in Interval Analysis?} 
\label{OtherNonTraditSect}

Next, let us turn to the other kinds of non-traditional intervals, such as improper 
(``reversed'') intervals, like $[2, 1]$, $[1, -2]$, etc. 
  
Some time ago, A.\,Neumaier reviewed in \cite{NeumaierDraft} the properties and 
applications of some non-traditional interval arithmetics (he called them the unfortunate 
term ``non-standard arithmetics'', as if someone issued a standard on various types 
of intervals). But since A.\,Neumaier himself is a pessimist and does not believe much 
in the usefulness of these arithmetics, then his review of applications turned out 
to be also pessimistic, almost like an obituary. Below, we try to give another overview 
of the capabilities of improper intervals from a more general point of view and show 
that, sooner or later, they will win their rightful place among the mathematical 
tools of natural and social science.

  
\subsection{Algebra I}

In the previous section, we considered the intervals from the viewpoint of the science 
of Topology. Let us now consider sets of intervals, more precisely, various interval 
arithmetics, from the viewpoint of another great science --- Algebra. 
  
The science of Algebra is often called the science of algebraic systems, i.\,e. 
the science that studies sets with certain operations and relations defined on them. 
Let us look at the operations existing on the set of intervals.
  
In terms of algebra, operations can be different. If an associative binary operation
is defined on the set of some elements, then this set is called a semigroup, a monoid,
or a group, depending on what properties of this operation are. Strictly speaking,
interval arithmetic is an algebraic system on which more than one operation is defined,
but for our analysis it is sufficient to consider these operations one at a time.
  
\emph{Semigroup} is the weakest formation, where almost nothing is required from 
a binary associative operation between elements. 
  
A \emph{monoid} is a semigroup with a neutral element. A reminder: neutral element 
or identity element is a special kind of element with respect to the binary operation 
on that set, which leaves other elements unchanged when combined with it.  
  
A \emph{group} is an algebraic system where the operation in question is reversible, 
that is, for any element, there is an inverse element with respect to this operation. 
  
In general, it is much more comfortable to work in a group than in a monoid or 
a semigroup. Implicit awareness of this fact has been one of the driving forces 
behind the expansion of popular and well-known algebraic systems over the past 
millennia. Recall that this is why the simplest natural numbers were once expanded 
to integers, then to rational numbers, and then to real and complex numbers, and 
so on (although this process was not linear). 
  
Why? The fact is, the operation in a group is ``predictable'' by its results and 
``invertible''. We can restore the operands from the result of the operation. We can 
perform algebraic manipulations in a group more easily and with less restrictions. 
In other words, our mathematical tools are richer in a group than in a semigroup 
or monoid. 
  
Specifically, in a group with the operation ``$\ast$'', if we have an equality 
\begin{equation*} 
a*c = b*c, 
\end{equation*} 
then we can conclude that 
\begin{equation*} 
a = b. 
\end{equation*} 
And if
\begin{equation*}
a*b = c, 
\end{equation*}
then
\begin{equation*}
a = b^{-1}*c, 
\end{equation*} 
where $b^{-1}$ is the inverse to $b$ with respect to the operation ``$\ast$''. 
Additionally, we can solve equations in the group, which is not possible 
in semigroups and monoids in the general case. 
  
Do we really need such capabilities in Interval Analysis? My answer is definitely 
``yes''. The author, for example, needs it, and he certainly knows that many others 
also need such things. The above is especially important when we solve the so-called 
``inverse problems'', when it is necessary to restore the preimage of a function 
by its value. A special case of ``inverse problems'' is the well-known problem 
of solving equations and systems of equations. 
  
If we cannot restore the preimage in elementary interval operations, then we do 
not have adequate tools to solve the ``inverse problems'' in general. 
  
Yet another obvious example where the above algebraic properties proves to be 
indispensable is metrology and measurement theory. This is the fundamental concept of 
\emph{measurement error}. Recall that by definition, an error is the difference between 
the approximate value of a quantity and its exact ideal value. In the natural sciences 
and engineering, this latter is understood as the true value of a physical quantity, 
that is, the value that ideally reflects the considered quantity or phenomenon within 
the framework of the model (theory) we have adopted to describe it. Anyway, 
the difference in the above formulation means an algebraic difference, i.e. addition 
with the opposite element with respect to addition. If the measurement result and/or 
the true value of a quantity are of an interval type, then it is not possible 
to correctly find the error in the classical interval arithmetic, since there is 
neither algebraic subtraction nor elements that are opposite to the proper intervals 
with respect to addition. 
  
In a similar situation in Geometry, when the Minkowski sum of sets is considered, and 
it is required to ``inverse'' it, the so-called Hukuhara difference was introduced 
\cite{Hukuhara} according to the following rule:
\begin{equation*} 
A\ominus B = C \qquad \Leftrightarrow \qquad A = B + C. 
\end{equation*} 
With respect to the classical interval arithmetic $\mbb{IR}$, it is better not to limit 
ourselves to partial reversal of operations, but to correct the situation fundamentally 
by performing algebraic completion of $\mbb{IR}$.

  
\subsection{Algebra II}

Let us consider further facts from Algebra. Even if an algebraic system with 
an associative binary operation is not a group, we can judge how good or bad 
this operation is in terms of its ``invertibility''. 
  
The condition
\begin{equation} 
\label{CanceLaw} 
a*c = b*c  \quad\Rightarrow \quad  a = b 
\end{equation} 
is called \emph{cancellation law}. If it holds in a semigroup or monoid, then this is
a sign that the operation ``$*$'' has good ``invertibility properties'', and it is 
almost as that from a group. Moreover, such a semigroup can often be enlarged to 
a group, or, in other words, this semigroup can be isomorphically embedded in a group. 

The corresponding result from Algebra is as follows: 
  
\bigskip\noindent  
\textbf{Theorem.} 
\textit{Every commutative semigroup that satisfies the cancellation law can be 
isomorphically embedded in an commutative group.} 
  
\bigskip 
A reader can see details, for example, in the book \cite{Kurosh}, Chapter 2, Section 5. 
The result we cite is at the bottom of page 47 of the Pergamon Press edition published 
in 1965 (its electronic version is available in Internet). 
   
Well, what about our interval arithmetic $\mbb{IR}$? 
   
It is an Abelian (commutative) semigroup, with respect to both addition and 
multiplication. For addition, the cancellation law is evidently satisfied,  but, 
for multiplication, the cancellation law is fulfilled only for intervals that do not 
contain zero. In the general case, the cancellation law is not valid, as one can see 
from the following example: 
\begin{equation*} 
[-1, 2]\cdot[2, 3] = [-3, 6] = [-1, 2]\cdot[1, 3]. 
\end{equation*} 
  
As a consequence, we can embed interval arithmetic in a broader algebraic system 
in which every element has an additive inverse (opposite) element, and any interval 
that does not contain zero, has the multiplicative inverse. 
 
This is the well-known Kaucher interval arithmetic $\mbb{KR}$, developed in PhD 
thesis of  Edgar  Kaucher \cite{Kaucher},  which  was  defended  in  Karlsruhe,  
Germany, in 1973, under the supervision of Prof.~Ulrich Kulisch. The main results 
of this dissertation were included in the articles \cite{Kaucher77, Kaucher80}. 
Earlier the idea of algebraic extension and completion of the classical interval 
arithmetic was also implemented in the preprint of H.-J. Ortolf \cite{Ortolf}, 
although it was not elaborated in detail. 
  
The operation opposite to the addition in $\mbb{KR}$, the so-called ``algebraic 
subtraction'', is an analog to the Hukuhara difference \cite{Hukuhara} and is usually 
denoted by the same symbol ``$\ominus$''. In the example of interval measurement error 
from the preceding subsection, we can therefore define it as the algebraic difference 
$(\tilde{\mbf{x}}\ominus\mbf{x}^\ast)$ between the measured value $\tilde{\mbf{x}}$ 
and the true value $\mbf{x}^\ast$ of a physical quantity. 
  
A lot is said about the Kaucher interval arithmetic in the standard IEEE 1788-2015 
for the implementation of interval arithmetic on digital computers \cite{IEEE1788}, 
although this arithmetic itself has yet to become a daily working tool for people 
using interval computation. 
  
It is necessary to say that E.\,Kaucher's work was quite nontrivial, since 
he had to extend interval multiplication with the help of not only algebraic 
considerations, but also with the use of the inclusion order relation, which was due 
to the lack of the multiplicative cancellation law. This also results in the fact 
that Kaucher interval arithmetic has some ``strange'' features, such as e.\,g. 
nontrivial zero divisors: 
\begin{equation*}
[-1, 1]\cdot[2, -3] = 0, 
\end{equation*}
which can be easily explained and interpreted from a more advanced standpoint. 
Namely, 
\begin{equation*}
[-1, 1]\cdot[2, -3] \ 
   = \  \bigvee_{x\in[-1, 1]} \  \bigwedge_{y\in[-3 ,2]} (x\cdot y) \  = \  0, 
\end{equation*} 
according to the min-max definition of arithmetic operations in the Kaucher interval 
arithmetic \cite{Kaucher, SharySurvey, SharyBook}. 
  
Moreover, as is often the case in mathematics, progress in one area immediately 
leads to advances in other areas related to the first one. For intervals in the new 
algebraically completed interval arithmetic, new logical interpretations of arithmetic 
operations are possible. They were discovered in the works of Spanish researchers 
in the 70-90s of the XX century and summarized in the book \cite{ModalIntAnal}. 
An alternative presentation of this theory can be found in the works 
\cite{Goldsztejn11, Goldsztejn22}. 
  
Additionally, the Kaucher arithmetic is also a lattice with respect to inclusion, 
but this is achieved in a more elegant way than simply assigning the empty set 
to the minimum of two nonintersecting intervals. Namely, in $\mbb{KR}$ 
\begin{equation*} 
\text{ minimum of $\mbf{a}$ and $\mbf{b}$ with respect to inclusion } \ 
   = \  \bigl[\,\max\{\un{\mbf{a}}, \un{\mbf{b}}\}, 
              \,\min\{\ov{\mbf{a}}, \ov{\mbf{b}}\}\,\bigr]. 
\end{equation*} 
If the intervals $\mbf{a}$ and $\mbf{b}$ do not intersect, the minimum is 
an ``improper interval''. 
  
Anyway, it makes sense to conclude this section by stressing the crucial role 
of the cancellation law in semigroups. As we have already said, this is a sign 
of partial "ivertibility" of the operation under study, and this fact greatly 
simplifies the solution of various inverse problems in specific semigroups.

  
\subsection{Algebra and beyond}

We can give the arguments of the previous subsections a slightly different context 
and show a different standpoint.
  
For several thousand years, there exists a very general and very powerful method 
for solving various mathematical problems, which is called the ``method of equations''. 
Its essence is 
\begin{list}{}{\itemsep=0pt\topsep=2pt\parsep=0pt} 
\item 
to designate the sought-for value through a special symbol \\ 
(usually a letter called ``unknown variable'') \\ 
and then 
\item 
to write out an equality (or several equalities, i.\,e., their system) \\ 
that the solution to the problem of interest must satisfy. 
\end{list} 
An equality with unknown variable whose value we have to find is called \emph{equation}. 
Further, to solve the original problem, it is necessary to solve the equation, i.\,e. 
find, in one way or another, the value of the unknown variable (which can be a number, 
a function, etc.) that satisfies the constructed equation or system of equations. 
  
The convenience and generality of this method is that the equation can be ``very 
implicit'' with respect to the unknown quantity. Moreover, the ways in which we search 
for its solution do not necessarily have to make meaningful practical sense with respect 
to the unknown variable. Instead, they can be very formal manipulations that are only 
mathematical in nature. It is only important that the resulting solution of the equation 
has a practical meaning. With the help of what kind of mathematics we got it, it is not 
so important. 
  
A nontrivial fact that some of readers (and even experts in Interval Analysis) may not 
realize: in Interval Analysis it is also useful to solve equations, Interval Equations. 
It is useful to find the solutions of interval equations in the general mathematical 
sense described above. We call them ``formal solutions'', since the nature of 
operations involved in the equation may be not necessarily algebraic. Anyway, doing 
this, of course, is best in an algebraically completed interval arithmetic, that is, 
in $\mbb{KR}$. 
  
``Formal solutions'' to interval equations were first considered in 1969, in the work 
of the Romanian mathematician S.\,Berti \cite{SBerti}, where they were not named in 
any way. Berti studied an interval quadratic equation and simply drew attention to the 
fact that the concept of solving an interval equation can also be given such a meaning. 
Then H.\,Ratschek and W.\,Sauer \cite{RatschekSauer} studied such solutions for a single 
interval linear equation, and they used the term ``algebraic solution''. In \cite{KNickel}, 
K.\,Nickel considered formal solutions to interval linear systems of equations in complex 
interval arithmetics, but did not name them in any specific way. Both the author himself 
and other researchers have also previously used the term ``algebraic solutions'' 
\cite{Markov-1999, Popova, Shary1996, SharyRC97}, but we now strongly recommend the term 
``formal solutions'' (see \cite{SharySurvey,ModalIntAnal,SharyBook,Shary-arXiv} and many 
others). 
  
For example, the interval $[0, 1]$ is a formal solution to the interval quadratic 
equation 
\begin{equation*} 
[1,2]\,x^2 + [-1,1]\,x = [-1,3]. 
\end{equation*} 
Interval function $\mbf{x}(t) = 10.5\,[\,e^t, e^{2t}\,]$ of the real argument $t$  
is a formal solution of the interval differential equation 
\begin{equation*}
\frac{dx(t)}{dt} = [1, 2]\;x(t).
\end{equation*}
The interval function $\mbf{x}(t) = [\,0, 2t\,]$ on $[0,1]$ is a formal solution 
to the Fredholm interval integral equation of the second kind 
\begin{equation*}
x(t) + \int_0^1 (1.5s+t)\,x(s)\, ds = [\,0, 3t+1\,]. 
\end{equation*}
The last two (purely illustrative) examples show the main drawback of the term 
``algebraic solution'': it emphasizes the algebraic nature of the operations 
that form the interval equation in question, so that talking about ``algebraic'' 
solution of interval differential, integral and such like equations is at least 
incorrect. 
  
Let us remind of some results obtained in 60-90's of the last century, showing 
the usefulness of ``formal solutions''. 
  
\paragraph{Enclosing the united solution set.} 
Let us be given an interval system of linear algebraic equations $\mbf{A}x = \mbf{b}$, 
with an interval $m\times n$-matrix $\mbf{A}$ and interval right-hand side $m$-vector 
$\mbf{b}$. Its united solution set is known to be the set 
\begin{equation*} 
\varXi_{\mathrm{uni}} (\mbf{A}, \mbf{b}) = 
   \bigl\{\,x\in\mbb{R}^n \mid Ax = b \; 
   \text{ for some $A\in\mbf{A}$ and $b\in\mbf{b}$}\;\bigr\} 
\end{equation*} 
i.\,e., the set of solutions to all point systems $Ax = b$ with $A\in\mbf{A}$ and 
$b\in\mbf{b}$. Interval estimation of the united solution set is an important 
practical problem, which is also one of the classic problems of Interval Analysis. 
Hundreds of articles have been devoted to it from the 60s of the last century 
up to the present. 
  
It is easy to show that the united solution set of the original system of equations 
coincides with the united solution set of the system in a fixed-point form 
\begin{equation*} 
x = (I - \mbf{A})\,x + \mbf{b}. 
\end{equation*} 
Next, a formal solution to the above fixed-point interval system gives an enclosure 
(outer interval box) of the united solution set, if the spectral radius of the matrix 
$|I - \mbf{A}|$, composed of the moduli of elements from $(I-\mbf{A})$, is less 
than $1$. This is the well-known result of Apostolatos and Kulisch 
\cite{ApostolatosKulisch}, obtained in 1968, which we reformulate in new terms 
convenient to our purposes. The reader can also see this result in the beginning 
of Chapter~12 of \cite{AlefeldHerzberger}. 
  
\paragraph{Inner estimation of the united solution set.} 
If an interval system of linear equations $\mbf{A}x = \mbf{b}$ is given, then 
a proper formal solution to the interval system 
\begin{equation*} 
(\dual\mbf{A})\,x = \mbf{b}, 
\end{equation*} 
where dual is dualization in Kaucher arithmetic, provides an inner box of the united 
solution set. This inner box is almost always inclusion maximal (that is, it touches 
the boundaries of the united solution set). 
  
\paragraph{Inner estimation of the tolerable solution set.} 
The tolerable solution set for an interval linear system $\mbf{A}x = \mbf{b}$ 
is known to be 
\begin{equation*} 
\varXi_{\mathrm{tol}} (\mbf{A}, \mbf{b}) = 
   \bigl\{\,x\in\mbb{R}^n \mid Ax\in\mbf{b} \; 
   \text{ for every $A\in\mbf{A}$}\;\bigr\} 
\end{equation*} 
i.\,e., the set of all such vectors $x$ that the product $Ax$ falls within 
the right-hand side vector $\mbf{b}$ for every $A\in\mbf{A}$. This is the second, 
in importance, among solution sets to interval systems of equations. 
  
Any proper formal solution to the interval system 
\begin{equation*} 
\mbf{A}x = \mbf{b} 
\end{equation*} 
(with the same form as the initial interval system) gives an inner box of the tolerable 
solution set. It is also inclusion maximal in most cases.

\paragraph{Enclosing the tolerable solution set.} 
Any proper formal solution to the interval system 
\begin{equation*} 
x = (I - (\dual\mbf{A}))\,x + \mbf{b} 
\end{equation*}
gives an enclosure of the tolerable solution set, if the spectral radius 
of $|I - \mbf{A}|$ is less than $1$. 
  
\bigskip 
And so on. The list is indeed very extensive, and we could continue it, but the above 
is enough for our short note. Naturally, there exist generalizations of the above 
results to nonlinear systems (see, e.\,g., \cite{ModalIntAnal}).

In conclusion, it is worth noting that, historically, the short term ``solution'' 
as applied to interval equations has taken on a slightly different meaning. Since 
the early 60s of the last century, when speaking about a \emph{solution} of an interval 
equation, one has been referring to the solution of some extended problem statement 
related to this equation. For example, ``to find an interval enclosure for the united 
solution set to an interval equation'' (a typical example of such terminology is the work 
\cite{PolyakNazin}). Or, ``to find an inner interval box within the tolerable solution 
set'' (this is the so-called interval tolerance problem). And so on. In other words, 
the situation at this point is similar to what we have in the theory of differential 
equations, where we do not talk about solutions to individual differential equations, 
per se. Usually some problem related to the differential equation in question is 
formulated (initial value problem, boundary value problem, etc.), which imposes 
additional conditions on the desired solution, without which the statement would be 
incomplete and meaningless. Then the solutions for this extended problem are considered, 
not for the single equation itself. The same is true for Interval Analysis.

  
\section{Empty intervals} 
\label{EmptyIntSect}

Empty intervals are useful in some cases, although this is mostly applicable 
in the classical interval arithmetic, equipped with set-theoretic operations. 
It usually results from the intersection, e.\,g.  
\begin{equation*} 
[1, 2]\cap[3, 4] = \varnothing. 
\end{equation*} 
In Kaucher interval arithmetic, it is sometimes advisable to use the minimum with 
respect to inclusion instead. Taking the minimum by inclusion is an operation similar 
in purpose to the intersection, but ``more friendly''. For example, 
\begin{equation*} 
[1, 2]\wedge[3, 4] = [3, 2], 
\end{equation*} 
and we thus get a non-empty result which can lead to nontrivial conclusions 
in the further reasoning. 
  
What happens when we add the empty set (``empty interval'') to an interval arithmetics? 
We have then 
\begin{equation} 
\label{EmptyInterval}
\begin{array}{ccccc} 
\mbf{a} + \varnothing &=& \varnothing + \mbf{a} &=& \varnothing,     \\[2mm] 
\mbf{a} - \varnothing &=& \varnothing - \mbf{a} &=& \varnothing,     \\[2mm] 
\mbf{a}\cdot\varnothing &=& \varnothing\cdot\mbf{a} &=& \varnothing, \\[2mm] 
\mbf{a} / \varnothing &=& \varnothing / \mbf{a} &=& \varnothing, 
\end{array} 
\end{equation}
and the cancellation law \eqref{CanceLaw} is ruined for both addition and multiplication 
in interval arithmetics. 
  
During the 2018 discussion, one of the participants, John Gustafson, compared the empty 
set to zero, i.\,e. to $0$. That was the wrong metaphor: unlike the noble identity 
elements like $0$ and $1$, the empty set is a kind of ``vampire'' in the algebraic 
sense, judging by equalities \eqref{EmptyInterval}. 
  
No invertibility of operations. No embedding in a larger and more complete 
algebraic system. The interval arithmetics with the empty set are suitable to solve 
mostly ``direct problems'' and perform chains of calculations in the forward direction. 
  
In particular, Kaucher interval arithmetic is incompatible with the empty set 
as we can see from the above reasons.

  
\section*{Acknowledgements}

The author thanks R. Baker Kearfott for his proposal to organize the outcomes 
of the on-line discussion to the present text.

\bibliographystyle{plain}
\bibliography{example}

  
\end{document}